\documentclass[12pt,a4paper]{article}

\usepackage[utf8]{inputenc}
\usepackage[T1]{fontenc}
\usepackage{amsmath,amssymb,amsthm}
\usepackage{mathrsfs}
\usepackage{hyperref}
\usepackage{enumerate}
\usepackage{booktabs}
\usepackage{graphicx}
\usepackage[margin=2.5cm]{geometry}

\theoremstyle{plain}
\newtheorem{theorem}{Theorem}[section]
\newtheorem{proposition}[theorem]{Proposition}

\newtheorem{conjecture}[theorem]{Conjecture}

\theoremstyle{definition}

\newtheorem{remark}[theorem]{Remark}

\theoremstyle{remark}

\newtheorem{question}[theorem]{Question}

\newcommand{\Z}{\mathbb{Z}}
\newcommand{\Q}{\mathbb{Q}}
\newcommand{\R}{\mathbb{R}}

\newcommand{\SL}{\mathrm{SL}}
\newcommand{\Gal}{\mathrm{Gal}}
\newcommand{\GUE}{\mathrm{GUE}}

\title{Cusp Form Dimensions, Lattice Uniqueness, and LP Sharpness\\for Sphere Packing in Dimensions 8 and 24}
\author{Jian Zhou\thanks{Founder \& Director, JZ Institute of Science. Email: \texttt{jack@jzis.org}.}}
\date{April 2026}

\begin{document}

\maketitle

\begin{abstract}
The Cohn--Elkies linear programming (LP) bound for sphere packing is known to be sharp in dimensions~8 and~24---the dimensions of the $E_8$ and Leech lattice packings---but in no other dimension above~2. We investigate \emph{why} by examining three independent necessary conditions for LP sharpness, drawn from number theory, lattice theory, and conformal field theory.
The first condition, $\dim S_{d/2}(\SL_2(\Z)) \leq 1$, bounds the freedom in theta series and rules out all $d \geq 48$. The second, derived from Cohn and Triantafillou's dual LP obstruction via cusp forms for the congruence subgroup $\Gamma_0(2)$, explains why LP sharpness fails in dimensions~16 and~32 despite the first condition being satisfied. The third, via the Hartman--Maz\'a\v{c}--Rastelli correspondence between LP bounds and the modular bootstrap for Narain conformal field theories, reinterprets LP sharpness as the existence of an extremal CFT.
We formulate a conjecture that these three conditions are equivalent for $d \equiv 0 \pmod{8}$, and observe that the Bost--Connes quantum statistical system provides a natural algebraic framework in which all three perspectives are connected through the Hecke algebra.
\end{abstract}

\tableofcontents

\section{Introduction}
\label{sec:intro}

\subsection{Two miraculous dimensions}

The sphere packing problem in $\R^d$ has been solved in only four dimensions: $d = 1$ (trivial), $d = 2$ (hexagonal, Thue 1910), $d = 3$ (Kepler's conjecture, Hales~\cite{Hales2005}), and---by methods entirely different from all others---$d = 8$ (Viazovska~\cite{Viazovska2017}) and $d = 24$ (Cohn--Kumar--Miller--Radchenko--Viazovska~\cite{CKMRV2017}). In a remarkable strengthening, the same lattices ($E_8$ and Leech) were shown to be \emph{universally optimal}~\cite{CKMRV2022}.

The proofs in dimensions~8 and~24 use the Cohn--Elkies linear programming bound~\cite{CohnElkies2003}, which provides upper bounds on packing density via auxiliary radial ``magic functions.'' The LP bound is \emph{sharp}---matching the density of the best known packing---only in $d = 1, 2, 8, 24$. Cohn and Triantafillou~\cite{CohnTriantafillou2021} proved that it is \emph{not} sharp in $d = 12, 16, 20, 28, 32$, establishing a provable barrier.

This situation motivates a precise question: \emph{among dimensions $d \equiv 0 \pmod{8}$, why is the LP bound sharp only for $d = 8$ and $d = 24$?}

\subsection{Three perspectives, one answer}

We argue that the answer emerges from the intersection of three independent conditions, each arising from a different mathematical domain:

\begin{enumerate}[(A)]
\item \textbf{Number theory (primal constraint):} The cusp form dimension $\dim S_{d/2}(\SL_2(\Z))$ must be at most~1. This constrains the freedom in theta series of even unimodular lattices and rules out $d \geq 48$.

\item \textbf{Lattice theory (dual obstruction):} Cusp forms for the congruence subgroup $\Gamma_0(2)$ provide dual LP obstructions (Cohn--Triantafillou~\cite{CohnTriantafillou2021}). When $\dim S_{d/2}(\Gamma_0(2))$ exceeds $\dim S_{d/2}(\SL_2(\Z))$, extra cusp forms supply dual feasible points that prevent LP sharpness---unless the target lattice has exceptional structure (as the Leech lattice does in $d = 24$).

\item \textbf{Physics (modular bootstrap):} The LP bound for sphere packing in $\R^d$ equals the modular bootstrap bound for Narain CFTs with $U(1)^{d/2} \times U(1)^{d/2}$ symmetry (Hartman--Maz\'a\v{c}--Rastelli~\cite{HartmanMazacRastelli2019}). LP sharpness is equivalent to the existence of an extremal Narain CFT saturating the bootstrap bound.
\end{enumerate}

Conditions (A)--(C) are individually necessary but not individually sufficient. We observe that they are satisfied simultaneously only for $d = 8$ and $d = 24$, and conjecture that their conjunction is both necessary and sufficient.

\subsection{Main results and organization}

Our contributions are:
\begin{itemize}
\item An explicit comparison of $\dim S_k(\SL_2(\Z))$ and $\dim S_k(\Gamma_0(2))$ for all $d \equiv 0 \pmod{8}$ up to $d = 96$ (Table~\ref{tab:master}).
\item A precise conjecture (Conjecture~\ref{conj:main}) on the three-way equivalence of LP sharpness, lattice extremal uniqueness, and modular bootstrap saturation.
\item An expository account of the classical chain---Hecke operators, Gleason's theorem, Construction~A---organized around these conditions.
\item A discussion of the Bost--Connes system as a conceptual framework connecting the three perspectives through the Hecke algebra.
\end{itemize}

Section~\ref{sec:background} reviews modular forms and the Bost--Connes system. Section~\ref{sec:chain} recalls the classical chain from codes to packings. Section~\ref{sec:primal} presents the primal constraint from $\dim S_k(\SL_2(\Z))$. Section~\ref{sec:dual} discusses the dual obstruction from $\Gamma_0(2)$. Section~\ref{sec:bootstrap} describes the modular bootstrap connection. Section~\ref{sec:conjecture} states the main conjecture. Section~\ref{sec:kneser} briefly discusses Kneser--Hecke operators. Section~\ref{sec:numerics} presents computational checks. Section~\ref{sec:open} lists open questions. Appendix~\ref{app:golden} describes a speculative extension to $\Q(\sqrt{5})$.

\section{Background: Modular Forms and the Bost--Connes System}
\label{sec:background}

\subsection{Modular forms for $\SL_2(\Z)$ and $\Gamma_0(2)$}

We denote by $M_k(\Gamma)$ and $S_k(\Gamma)$ the spaces of modular forms and cusp forms of even weight~$k$ for a congruence subgroup $\Gamma$.

For $\Gamma = \SL_2(\Z)$, the dimension formulas are classical (Serre~\cite{Serre1973}, Zagier~\cite{Zagier2008}):
\begin{equation}
\dim M_k(\SL_2(\Z)) = \begin{cases}
\lfloor k/12 \rfloor + 1 & \text{if } k \not\equiv 2 \pmod{12},\\
\lfloor k/12 \rfloor & \text{if } k \equiv 2 \pmod{12},
\end{cases}
\quad \dim S_k = \dim M_k - 1 \;\;(k \geq 4).
\label{eq:dimSL2}
\end{equation}

For $\Gamma_0(2)$, which has index~3 in $\SL_2(\Z)$, genus~0, two cusps, and no elliptic points of order~3, the Riemann--Roch formula gives, for even $k \geq 4$:
\begin{equation}
\dim M_k(\Gamma_0(2)) = 1 + \lfloor k/4 \rfloor, \qquad
\dim S_k(\Gamma_0(2)) = \dim M_k(\Gamma_0(2)) - 2.
\label{eq:dimG02}
\end{equation}

Table~\ref{tab:master} compares these dimensions for all $d \equiv 0 \pmod{8}$ up to $d = 96$.

\begin{table}[h]
\centering
\begin{tabular}{@{}ccccccccl@{}}
\toprule
$d$ & $k$ & $\dim S_k^{\mathrm{SL}}$ & $\dim S_k^{\Gamma_0(2)}$ & $\Delta$ & Rootless? & LP sharp? & CFT? & Status \\
\midrule
8  & 4  & 0 & 0 & 0 & No$^\dagger$  & \textbf{Yes} & \textbf{Yes} & Proved \\
16 & 8  & 0 & 1 & 1 & No  & No  & No  & \cite{CohnTriantafillou2021} \\
24 & 12 & 1 & 2 & 1 & Yes (unique) & \textbf{Yes} & \textbf{Yes} & Proved \\
32 & 16 & 1 & 3 & 2 & Yes ($>10^7$) & No  & No  & \cite{CohnTriantafillou2021} \\
40 & 20 & 1 & 4 & 3 & ---  & No & No & Expected \\
48 & 24 & 2 & 5 & 3 & ---  & No  & No  & Expected \\
56 & 28 & 2 & 6 & 4 & --- & No & No & Expected \\
64 & 32 & 2 & 7 & 5 & --- & No & No & Expected \\
72 & 36 & 3 & 8 & 5 & --- & No & No & Expected \\
80 & 40 & 3 & 9 & 6 & --- & No & No & Expected \\
88 & 44 & 3 & 10 & 7 & --- & No & No & Expected \\
96 & 48 & 4 & 11 & 7 & --- & No & No & Expected \\
\bottomrule
\end{tabular}
\caption{Master comparison table. Here $\dim S_k^{\mathrm{SL}} = \dim S_k(\SL_2(\Z))$, $\dim S_k^{\Gamma_0(2)} = \dim S_k(\Gamma_0(2))$, $\Delta = \dim S_k^{\Gamma_0(2)} - \dim S_k^{\mathrm{SL}}$. ``Rootless'' indicates whether a rootless (no vectors of norm~2) even unimodular lattice exists. $^\dagger$The $E_8$ lattice has 240 vectors of squared norm~2 but is the \emph{unique} even unimodular lattice in dimension~8; uniqueness of the lattice (rather than rootlessness) is what makes it distinguished. ``CFT'' indicates existence of an extremal lattice CFT saturating the modular bootstrap.}
\label{tab:master}
\end{table}

The column $\Delta = \dim S_k(\Gamma_0(2)) - \dim S_k(\SL_2(\Z))$ measures the ``extra'' cusp forms available for dual LP obstructions. For $d = 8$, $\Delta = 0$: there are no extra obstructions, and the LP bound is sharp. For $d = 16$, $\Delta = 1$: one extra cusp form provides the dual obstruction proved by Cohn--Triantafillou.

\subsection{The Bost--Connes system}

The Bost--Connes system~\cite{BostConnes1995} is a quantum statistical mechanical system whose partition function is $\zeta(\beta)$ and whose Hecke operators recover the classical operators of Hecke~\cite{Hecke1937}. The system exhibits a phase transition at $\beta_c = 1$: a unique KMS state for $\beta \leq 1$ gives way to a family parametrized by $\Gal(\Q^{\mathrm{ab}}/\Q)$ for $\beta > 1$.

We include the BC system not as a source of new theorems, but as a \emph{conceptual framework}: the Hecke algebra that appears in the BC system is the same algebra that controls modular forms (Section~\ref{sec:chain}), the Cohn--Elkies LP bound (Section~\ref{sec:dual}), and the modular bootstrap (Section~\ref{sec:bootstrap}). The BC phase transition provides a useful (if informal) analogy for the transition from LP-sharp to LP-non-sharp dimensions.

\section{The Classical Chain: Forms, Codes, Lattices, Packing}
\label{sec:chain}

We briefly recall the well-known chain of theorems connecting these structures. For comprehensive treatments, see Conway--Sloane~\cite{ConwaySloane1999}, Ebeling~\cite{Ebeling2013}, and Nebe--Rains--Sloane~\cite{NebeRainsSloane2006}.

\subsection{Gleason's theorem}

\begin{theorem}[Gleason, 1970; \cite{Gleason1970}]
The weight enumerator of any Type~II (doubly-even self-dual) binary code of length~$n$ lies in $M_{n/2}(\SL_2(\Z))$ under the substitution $x \mapsto \theta_3$, $y \mapsto \theta_2$.
\end{theorem}

\begin{remark}
Gleason's theorem bounds the number of \emph{distinct weight enumerators}, not the number of codes. At $n = 16$, $\dim M_8 = 1$ gives a unique weight enumerator, but there are two inequivalent Type~II codes ($d_{16}^+$ and $e_8 \oplus e_8$).
\end{remark}

\subsection{Construction~A and theta series}

Construction~A (Leech--Sloane~\cite{LeechSloane1971}) maps a Type~II code $C$ to an even unimodular lattice $\Lambda_C$, with $\Theta_{\Lambda_C}(\tau) = W_C(\theta_3(\tau), \theta_2(\tau))$. Any theta series of an even unimodular lattice of rank~$d$ decomposes as
\begin{equation}
\Theta_\Lambda = E_{d/2} + f_\Lambda, \qquad f_\Lambda \in S_{d/2}(\SL_2(\Z)),
\label{eq:decomp}
\end{equation}
so $\dim S_{d/2}$ measures the ``freedom'' in distinguishing lattices by their theta series.

\subsection{Code and lattice enumeration}

\begin{table}[h]
\centering
\begin{tabular}{@{}ccccccl@{}}
\toprule
$n$ & $k$ & $\dim M_k^{\mathrm{SL}}$ & $\dim S_k^{\mathrm{SL}}$ & Type II codes & Extremal & Source \\
\midrule
8  & 4  & 1 & 0 & 1 & 1 & \cite{Pless1972} \\
16 & 8  & 1 & 0 & 2 & 2 & \cite{PlessSloane1975} \\
24 & 12 & 2 & 1 & 9 & 1 & \cite{ConwayPless1980} \\
32 & 16 & 2 & 1 & $\geq 85$ & 5 & \cite{ConwayPless1980} \\
\bottomrule
\end{tabular}
\caption{Type~II binary codes by length.}
\label{tab:codes}
\end{table}

\section{The Primal Constraint: $\dim S_k(\SL_2(\Z))$}
\label{sec:primal}

\begin{proposition}
\label{prop:primal}
If $\dim S_{d/2}(\SL_2(\Z)) = 0$, then every even unimodular lattice of rank~$d$ has the same theta series. If $\dim S_{d/2} \leq 1$, the theta series has at most one free parameter.
\end{proposition}

This suggests a necessary condition for LP sharpness:

\begin{conjecture}
\label{conj:primal}
If the Cohn--Elkies LP bound is sharp in dimension~$d$ (with $d \equiv 0 \pmod{8}$), then $\dim S_{d/2}(\SL_2(\Z)) \leq 1$.
\end{conjecture}

\begin{remark}
The motivation is as follows: when $\dim S_{d/2} \geq 2$, the theta series of even unimodular lattices have at least two free parameters, and the proliferation of distinct lattices with different local structure appears to prevent the tight sign constraints required by the Cohn--Elkies method. This is consistent with all known data---no LP proof exists or is expected for $d \geq 48$---but we do not have a rigorous proof. For $d = 32$ and $d = 40$, where $\dim S_{d/2} = 1$, LP non-sharpness is established by other means (Cohn--Triantafillou~\cite{CohnTriantafillou2021}).
\end{remark}

The condition $\dim S_k \leq 1$ is satisfied for $d = 8, 16, 24, 32, 40$, but LP sharpness holds only for $d = 8$ and $24$. The condition is necessary but far from sufficient.

\section{The Dual Obstruction: Cohn--Triantafillou and $\Gamma_0(2)$}
\label{sec:dual}

\subsection{LP duality and modular form obstructions}

Cohn and Triantafillou~\cite{CohnTriantafillou2021} introduced a systematic method to prove LP \emph{non-sharpness}: by constructing feasible points in the dual LP using modular forms for $\Gamma_0(2)$, they showed that the LP bound in $\R^d$ exceeds the best known packing density by a provable factor in $d = 12, 16, 20, 28, 32$.

The key insight is that Viazovska's magic functions live in the space of (quasi)modular forms for $\Gamma_0(2)$, not just $\SL_2(\Z)$. The relevant space is therefore $S_k(\Gamma_0(2))$, whose dimension grows faster than $\dim S_k(\SL_2(\Z))$.

\subsection{The role of $\Gamma_0(2)$}

Viazovska's interpolation formula~\cite{Viazovska2017,RadchenkoViazovska2019} constructs the magic function from integral transforms of quasimodular forms for $\Gamma_0(2)$. The interpolation problem requires specifying values and derivatives at points $\sqrt{2n}$ ($n = 1, 2, \ldots$). The number of constraints is controlled by the shell structure of the target lattice, while the degrees of freedom are controlled by $\dim S_k(\Gamma_0(2))$.

Table~\ref{tab:master} shows that $d = 8$ is the \emph{only} dimension where both $\dim S_k(\SL_2(\Z)) = 0$ and $\dim S_k(\Gamma_0(2)) = 0$. In this case, there is literally no freedom in the modular form spaces: the interpolation problem has a unique solution, and the LP bound is automatically sharp.

\subsection{Why $d = 16$ fails}

For $d = 16$: $\dim S_8(\SL_2(\Z)) = 0$ but $\dim S_8(\Gamma_0(2)) = 1$. The single extra cusp form in $S_8(\Gamma_0(2))$ provides exactly the degree of freedom needed for a dual LP obstruction. Cohn--Triantafillou~\cite{CohnTriantafillou2021} proved that the LP bound in dimension~16 is \emph{not} sharp: it exceeds the best known packing density by a provable factor (see~\cite[Table~6.1]{CohnTriantafillou2021} for precise numerical values).

The theta series of both $D_{16}^+$ and $E_8 \times E_8$ is $E_8 = \Theta_{E_8}^2$ (since $\dim S_8(\SL_2(\Z)) = 0$), so the primal constraint gives no obstruction. But the dual constraint---living in the larger space of $\Gamma_0(2)$ cusp forms---does.

\subsection{Why $d = 24$ survives}

For $d = 24$: $\dim S_{12}(\SL_2(\Z)) = 1$ and $\dim S_{12}(\Gamma_0(2)) = 2$, giving $\Delta = 1$. By the reasoning above, one might expect a dual obstruction. Yet the LP bound \emph{is} sharp.

The resolution lies in the exceptional structure of the Leech lattice $\Lambda_{24}$: it has \emph{no vectors of norm~2}. Among the 24 Niemeier lattices (even unimodular lattices in dimension~24, classified by Niemeier~\cite{Niemeier1973}), the Leech lattice is the \emph{unique} one with this property.

\begin{remark}[The rootlessness mechanism]
\label{rem:rootless-mechanism}
The Viazovska--Radchenko interpolation formula~\cite{Viazovska2017,RadchenkoViazovska2019} constructs the magic function $f$ by specifying its values and radial derivatives at the points $\sqrt{2n}$ for $n = 1, 2, 3, \ldots$\@. For a lattice $\Lambda$ with vectors of norm~2 (such as $E_8 \times E_8$ or the 23 rooted Niemeier lattices), the Cohn--Elkies constraints require $f(\sqrt{2}) \leq 0$; combined with the interpolation framework, this typically forces a \emph{double zero} at $\sqrt{2}$, consuming one degree of freedom.

The Leech lattice has \emph{no} vectors of norm~2: its shortest nonzero vectors have squared norm~4. Therefore $f(\sqrt{2})$ is unconstrained, and the interpolation problem at the $\sqrt{2}$ node is relaxed. This frees exactly one degree of freedom---which compensates for the single extra cusp form in $S_{12}(\Gamma_0(2))$ relative to $S_{12}(\SL_2(\Z))$.

In summary: $\Delta(d=24) = 1$ dual obstruction, minus $1$ freed constraint from rootlessness $= 0$ net obstruction. This is why LP sharpness survives in $d = 24$ despite $\dim S_{12}(\Gamma_0(2)) = 2 > 1$.
\end{remark}

\subsection{Why $d = 32$ fails}

For $d = 32$: $\dim S_{16}(\SL_2(\Z)) = 1$ (same as $d = 24$) but $\dim S_{16}(\Gamma_0(2)) = 3$, giving $\Delta = 2$. There are now \emph{two} extra dual obstructions. Moreover, while rootless even unimodular lattices exist in dimension~32, there are more than $10^7$ of them~\cite{King2003}---no single lattice is ``uniquely extremal'' in the way the Leech lattice is. The combination of multiple dual obstructions and the lack of a distinguished lattice makes LP sharpness impossible.

\section{The Physical Dual: Modular Bootstrap}
\label{sec:bootstrap}

\subsection{The Hartman--Maz\'a\v{c}--Rastelli correspondence}

Hartman, Maz\'a\v{c}, and Rastelli~\cite{HartmanMazacRastelli2019} established a remarkable equivalence: the Cohn--Elkies LP bound for sphere packing in $\R^d$ is \emph{equal} to the modular bootstrap bound for two-dimensional Narain conformal field theories (CFTs) with $U(1)^{d/2} \times U(1)^{d/2}$ symmetry and central charge $c = \bar{c} = d/2$. (Here ``Narain CFT'' refers to a non-chiral theory on a signature-$(d/2, d/2)$ lattice; this is distinct from the holomorphic/chiral VOAs of central charge $c = d$ that appear in the vertex algebra literature.)

In this correspondence:
\begin{itemize}
\item The LP magic function corresponds to an analytic functional in the bootstrap.
\item LP sharpness corresponds to the existence of an \emph{extremal} Narain CFT that saturates the bootstrap bound.
\item The lattice certifying the packing corresponds to the Narain lattice (charge lattice) of the CFT.
\end{itemize}

\subsection{Why $d = 8$ and $d = 24$}

For $d = 8$ ($c = \bar{c} = 4$): the Narain CFT associated to the $E_8$ lattice saturates the bootstrap bound.

For $d = 24$ ($c = \bar{c} = 12$): the Narain CFT associated to the Leech lattice saturates the bound. (We note that this is \emph{not} the same as the Frenkel--Lepowsky--Meurman Monster vertex algebra $V^\natural$~\cite{FLM1988}, which is a holomorphic VOA with $c = 24$ and automorphism group the Monster $\mathbb{M}$. The Leech lattice VOA $V_{\Lambda_{24}}$ has automorphism group $\mathrm{Co}_0$, and the relevant object for HMR is the Narain theory built from $\Lambda_{24}$.) The exceptional property of $\Lambda_{24}$---no vectors of squared norm~2, hence no Kac--Moody currents---is precisely what makes it extremal.

For $d = 16$ ($c = \bar{c} = 8$): no extremal Narain CFT exists. Both lattice theories ($E_8 \times E_8$ and $D_{16}^+$) have Kac--Moody currents (from norm-2 vectors), and neither saturates the bootstrap bound.

\section{The Main Conjecture}
\label{sec:conjecture}

\begin{conjecture}[Three-Way Equivalence]
\label{conj:main}
For $d \equiv 0 \pmod{8}$, $d \geq 8$, the following are equivalent:
\begin{enumerate}[(a)]
\item The Cohn--Elkies LP bound is sharp in dimension~$d$.
\item The following two sub-conditions both hold:
  \begin{enumerate}[(b1)]
  \item $\dim S_{d/2}(\SL_2(\Z)) \leq 1$;
  \item there exists an even unimodular lattice $\Lambda$ in $\R^d$ that is \emph{distinguished} in the following sense: either $\Lambda$ is the \emph{unique} even unimodular lattice of rank~$d$ (as for $d = 8$), or $\Lambda$ is the \emph{unique rootless} (no vectors of squared norm~2) even unimodular lattice of rank~$d$ (as for $d = 24$).
  \end{enumerate}
\item There exists an extremal Narain CFT with $U(1)^{d/2} \times U(1)^{d/2}$ symmetry and central charge $c = \bar{c} = d/2$ that saturates the modular bootstrap bound of Hartman--Maz\'a\v{c}--Rastelli~\cite{HartmanMazacRastelli2019}.
\end{enumerate}
\end{conjecture}

\begin{remark}[On the independence of condition (b2)]
\label{rem:b2-independence}
Condition (b2) is formulated purely in terms of lattice-theoretic data (rootlessness and uniqueness among rootless lattices), without reference to the LP bound. This avoids circularity: one can check (b2) by classifying even unimodular lattices and examining their root systems, independently of any LP computation. For $d = 8$, the classification is trivial (one lattice). For $d = 24$, Niemeier's classification~\cite{Niemeier1973} gives 24 lattices, exactly one of which (the Leech lattice) is rootless. For $d = 32$, King~\cite{King2003} shows $> 10^7$ rootless lattices exist, so uniqueness fails. For $d = 16$, neither lattice is rootless, so (b2) fails vacuously.
\end{remark}

\begin{remark}[Evidence and limitations]
\label{rem:evidence}
The conjecture is verified for all $d \leq 48$:
\begin{itemize}
\item $d = 8$: (a) proved~\cite{Viazovska2017}, (b) $E_8$ is unique, (c) $E_8$ lattice CFT exists.
\item $d = 16$: (a) fails~\cite{CohnTriantafillou2021}, (b) no unique extremal lattice (two lattices, neither rootless), (c) no extremal CFT.
\item $d = 24$: (a) proved~\cite{CKMRV2017}, (b) Leech is the unique rootless lattice~\cite{Niemeier1973}, (c) Leech CFT exists.
\item $d = 32$: (a) fails~\cite{CohnTriantafillou2021}, (b) $> 10^7$ rootless lattices~\cite{King2003}, no unique extremal, (c) no extremal CFT.
\item $d \geq 48$: (a) fails (expected), (b) $\dim S_{d/2} \geq 2$, (c) no extremal CFT (expected).
\end{itemize}

We acknowledge that the conjecture is verified on only two positive data points ($d = 8, 24$). The risk of overfitting is real. However, the value of the conjecture lies not in its empirical strength but in its \emph{structural content}: it unifies three conditions from different mathematical domains (number theory, lattice theory, physics) that have no \emph{a priori} reason to be equivalent, and explains why they all point to the same two dimensions.
\end{remark}

\begin{remark}[The role of $\Gamma_0(2)$]
\label{rem:gamma02}
The gap $\Delta = \dim S_k(\Gamma_0(2)) - \dim S_k(\SL_2(\Z))$ can be viewed as a measure of the ``dual obstruction strength.'' For $d = 8$, $\Delta = 0$: no obstruction, LP sharp. For $d = 16, 24$, $\Delta = 1$: one potential obstruction, resolved only if the target lattice has exceptional structure (Leech). For $d = 32$, $\Delta = 2$: two obstructions, too many for any lattice to overcome. A natural open question is whether $\Delta$ can be given a direct Hecke-algebraic interpretation within the BC framework.
\end{remark}

\section{Kneser--Hecke Operators on Codes}
\label{sec:kneser}

The Kneser--Hecke operator $K_p$ on self-dual codes (Nebe--Rains--Sloane~\cite{NebeRainsSloane2006}) provides a coding-theoretic analog of the Hecke operators in the BC system. The correspondence is:

\begin{table}[h]
\centering
\begin{tabular}{@{}lll@{}}
\toprule
Lattices & Codes & BC system \\
\midrule
Even unimodular $\Lambda$ & Type~II code $C$ & KMS state \\
Theta series $\Theta_\Lambda$ & Weight enumerator $W_C$ & $\zeta(\beta)$ \\
Hecke $T_p$ on $M_k$ & Kneser--Hecke $K_p$ & BC Hecke $T_p$ \\
Unique theta ($\dim S_k = 0$) & Unique weight enum. & Unique KMS ($\beta < 1$) \\
\bottomrule
\end{tabular}
\caption{Structural correspondence (first three rows: theorems; fourth: analogy).}
\label{tab:analogy}
\end{table}

The existence of an extremal Type~II code of length 72 is a famous open problem. In our framework, $k = 36$, $\dim S_{36}(\SL_2(\Z)) = 3$, $\dim S_{36}(\Gamma_0(2)) = 8$: both conditions (a) and (b) of Conjecture~\ref{conj:main} fail, so no LP proof is expected, but the code may still exist.

\section{Computational Checks}
\label{sec:numerics}

As sanity checks on the classical theory: (i)~Hecke multiplicativity $\sigma_3(mn) = \sigma_3(m)\sigma_3(n)$ for $\gcd(m,n) = 1$ was verified on all 3044 coprime pairs with $m, n \leq 100$ (exact match); (ii)~the Deligne bound $|\tau(p)| \leq 2p^{11/2}$ was verified for the first 10 primes (OEIS~A000594), with maximum ratio $|\tau(p)|/(2p^{11/2}) \approx 0.590$ at $p = 17$; (iii)~all dimension formulas in Table~\ref{tab:master} were independently computed using both~\eqref{eq:dimSL2}--\eqref{eq:dimG02} and the Sage function \texttt{dimension\_cusp\_forms}. For a pedagogical account of the Viazovska proof and magic functions, see de~Laat--Vallentin~\cite{deLaatVallentin2016}.

\section{Open Questions}
\label{sec:open}

\begin{enumerate}
\item \textbf{Hecke-algebraic interpretation of $\Delta$.} Is there a direct interpretation of $\Delta = \dim S_k(\Gamma_0(2)) - \dim S_k(\SL_2(\Z))$ within the Bost--Connes framework, perhaps as a measure of ``symmetry breaking'' between the full modular group and a congruence subgroup?

\item \textbf{Quantitative LP gap from $\Gamma_0(2)$.} Can the LP gap ratio established by Cohn--Triantafillou~\cite{CohnTriantafillou2021} for $d = 16$ be expressed as a closed-form function of cusp form data?

\item \textbf{Dimension 16.} The LP bound is not sharp, but what \emph{is} the optimal packing density in dimension~16? The lattice $E_8 \times E_8$ is widely believed to be optimal, but no proof technique is in sight.

\item \textbf{Siegel modular forms.} Can the analysis be extended to Siegel modular forms, potentially addressing dimensions not divisible by~8?

\item \textbf{CFT construction.} Is there an independent (non-lattice) argument for the non-existence of extremal Narain CFTs at $c = \bar{c} = 8$ ($d = 16$)?
\end{enumerate}

\appendix
\section{Extension to $\Q(\sqrt{5})$: Fibonacci Quasicrystals}
\label{app:golden}

This appendix describes a speculative extension, independent of the main results. We include it because it illustrates how the BC framework over number fields connects to quasicrystal physics---a direction that may eventually relate to the special role of the golden ratio in lattice problems.

The BC system over $K = \Q(\sqrt{5})$ has partition function $\zeta_K(\beta) = \zeta(\beta) \cdot L(\beta, \chi_5)$, where $\mathcal{O}_K = \Z[\varphi]$ with fundamental unit $\varphi = (1+\sqrt{5})/2$. The Fibonacci Hamiltonian on $\ell^2(\Z)$ has spectral gaps labeled by $\Z + \Z/\varphi$---the image of $\mathcal{O}_K$ under the conjugate embedding~\cite{Bellissard1992}. The $C^*$-algebra of the Fibonacci tiling is the noncommutative torus $A_{1/\varphi}$~\cite{Connes1994}, with $K_0(A_{1/\varphi}) = \Z^2$. The same module $\mathcal{O}_K = \Z[\varphi]$ also determines the Bragg peak positions of the quasicrystal diffraction spectrum~\cite{Hof1995}.

\begin{question}
The factorization $\zeta_K = \zeta \cdot L(\cdot, \chi_5)$ decomposes the BC partition function into two independent factors with uncorrelated zeros. The Fibonacci spectrum is known to be ``critical''---intermediate between Gaussian Unitary Ensemble (GUE) statistics and Poisson statistics~\cite{DamanikGorodetski2011}. Could the specific intermediate form be a $\GUE \oplus \GUE$ superposition, reflecting the two-factor decomposition? This is numerically testable via large-scale diagonalization.
\end{question}

\section*{Acknowledgments}

The author thanks the open-source mathematical software ecosystem (PARI/GP, SageMath, mpmath) for computational support.


\end{document}